\documentclass[12pt]{article}
\usepackage{amsmath,amsthm,amsfonts,amssymb,amscd}
\usepackage[latin1]{inputenc}
\usepackage[brazil]{babel}

\newtheorem{theorem}{Theorem}[section]

\newtheorem{lemma}{Lemma}

\newcommand{\ddd}{\displaystyle}
\newcommand{\nd}{\noindent}
\newcommand{\vsp}{\vspace{0.2cm}}
\newcommand{\fimaf}{$\hfill{\Box}$}
\newcommand{\fim}{$\hfill{\rule{2.5mm}{2.5mm}}$}
\begin{document}

\title{\centering  On the study of a class of non-linear differential equations on compact Riemannian Manifolds.}

\maketitle

\begin{center}{\bf Carlos Silva, Romildo Pina \\ \bf and  Marcelo Souza }
 
\end{center}

\begin{abstract} 
{ We study the existence of solutions of the non-linear differential equations on the compact Riemannian manifolds $(M^n,g),\,\,
 n\geq 2$, 
 \begin{equation}
\label{E2}
\Delta_p u + a(x)u^{p-1} = \lambda f(u,x),
\end{equation} 
where $\Delta_p$ is the $p-$laplacian, with $1<p<n$.
 The equation (\ref{E2}) generalizes a equation considered by Aubin \cite{aubin3}, where he has considered,
  a compact Riemannian manifold  \ $(M,g)$,\  the  differential equation   ($p=2$) 
\begin{equation}
\label{E1}
\Delta u \ +\  a(x)u \ = \ \lambda f(u,x),
\end{equation}
where \ $a(x)$\  is a \ $C^{\infty}$ function  defined on \ $M$\  and \ $f(u,x)$\  is a \ $C^{\infty}$ function  defined on 
\ $\mathbb{R}\times M$.                          
We show that the equation \ (\ref{E2}) \ has solution \ $(\lambda,u)$, 
\ where \ $\lambda \ \in 
\ \mathbb{R}$, \ $u \ \geq \ 0$, \ $u \ \not\equiv \ 0$ \ is a function \ $C^{1,\alpha}$, \ $0 < \ \alpha \ < \  1$,
if \ $f \in C^{\infty}$ \ satisfies some
growth and parity conditions.}
\end{abstract}

keywords {critical Sobolev exponent, compact Riemannian manifold, non-linear differential equation}

Subjclass {Primary 53C21; Secondary 35J60}

{The authors were partially supported by CAPES}

\section{Introduction}

The study of the theory of nonlinear differential equations on 
Riemannian manifolds has began in 1960 with the so-called Yamabe problem. At a time when little was known 
about the methods of studying a non-linear equation, the Yamabe problem came to light of a geometric idea and from time sealed a 
merger of the areas of geometry and differential equations.
Let \ $(M,g)$ be a compact Riemannian manifold of dimension \ $n$, \ $n \geq 3$. \ Given \ $\widetilde{g} = u^{4/(n - 2)}g$
 \  some conformal metrical to the metric \ $g$, \ is well known that the scalar curvatures
  \ $R$ \ and \ $\widetilde{R}$ \ of the metrics \ $g$ \ and \ $\widetilde{g}$, 
respectively, satisfy the law of transformation

$$\Delta u \ + \ \ddd\frac{n - 2}{4(n - 1)}Ru \ = \ \ddd\frac{n - 2}{4(n-1)}\widetilde{R}u^{2^{\ast}-1}$$

\nd where \ $\Delta$ \ denote the Laplacian operator associated to \ $g$.

In 1960, Yamabe \cite{yamabe} announced that for every compact Riemannian manifold \ $(M,g)$ 
\ there exist  a metric \ $\widetilde{g}$ \ conformal to \ $g$ \ for which \ $\widetilde{R}$ \ is constant.
 In another words, this mean that for every compact Riemannian manifold \ $(M,g)$ there exist 
\ $u \in C^{\infty}(M), \ u > 0 \ \mbox{on} \ M$ \ and \ $\lambda \in \mathbb{R}$ \ such that  
\vspace{-0.5cm}

$$\Delta u \ + \ \ddd\frac{n - 2}{4(n - 1)}Ru \ = \ \lambda u^{2^{\ast}-1}. \eqno {(Y)}$$

In 1968, Trüdinger \cite{trudinger} found an error in the work of Yamabe, which generated a race to
solve what became known as the Yamabe problem, today it is completely positively solved,
that is, the assertion of  Yamabe is true.

The main step towards the resolution of the Yamabe problem was given in 1976 by T. Aubin in his classic article \cite{aubin2}.
 In \cite{aubin2} Aubin showed that the statement was true since the manifold satisfy a condition
on an invariant (called Yamabe invariant).
Then he used  tests functions, locally defined, to show that  non locally conformal flat manifolds, of
dimension  $ n \geq $ 6, \ satisfied this condition. Finally, for $n\geq 3 $ the problem was completed solved by 
R. Schoen \cite{schoen}. 

As previously reported, several disturbances were considered to the Yamabe's problem, all disturbances of analytical character, 
both in the sense of equation (with the addition of other factors) and in the sense of the operator (the Laplacian for
 the $ p $-Laplacian), and using the Aubin's idea of estimating corresponding functional. We
 can cite some articles, such as \cite {azorero1}, \cite {brezis1}, \cite {demegel}, \cite {djadli}, \cite {druet4} and 
\cite {olimpio}.

In \cite{silva2}, the author studied the existence of solutions for a class of non-linear differential  equation 
on compact Riemannian manifolds. 
 He establish a lower and upper solutions' method to show the existence of a smooth positive solution for the equation (\ref{E4})
\begin{equation}
 \label{E4}
\Delta u \ + \ a(x)u \ = \ f(x)F(u) \ + \ h(x)H(u),
\end{equation}

\nd where \ $a, \ f, \ h$ \  are positive smooth functions on $M^n$, a $n-$dimensional compact Riemannian manifold,
 and \ $ F, \ H$ 
\ are non-decreasing smooth functions on $\mathbb{R}$.  In \cite{djadli} the equation (\ref{E4}) was studied when
$F(u)=u^{2^{\ast}-1} $ and $H(u)=u^q$ in the Riemannian context, i.e.,
\begin{equation}
\label{E3}
\Delta u \ + \ a(x)u \ = \ f(x)u^{2^{\ast}-1} \ + \ h(x)u^q, 
\end{equation}
\nd where \ $0 \ < \ q \ < 1$. In \cite{correa} Corr\^ea,  Gon\c{c}alves and Melo  studied
an equation of the type equation (\ref{E3}),  in the Euclidean context,  with respect to a more general operator than the 
laplacian operator.

This work, which is organized into four sections, also aims to work with problems related to the equation 
$ (Y) $, although, as we shall see, with different methods from those used by Yamabe, these results were obtained in \cite{silva},

In section 2, we enter what we consider as basic concepts necessary to understand it, as some definitions and theorems of embedded.

\noindent We consider
$F(t,x) \ = \ \ddd\int^{t}_0{f(s,x)ds},$  
$B(u)\ = \ \ddd\int_M{F(u(x),x)dV}$ 

\noindent and

$I(u)\ =\ \ddd\int_M{|\nabla u|^pdV}\ +\ \ddd\int_M{a|u|^pdV}.$

\noindent Given \ $R > 0$, \  we also consider 
$\mathcal{H} \ = \ \{u \ \in \ H^p_1(M);\  B(u)\ = \ R \}$
  and 
$\mu_R \ = \ \ddd\inf_{u \in \mathcal{H}}{I(u)}.$

We proved, in the following theorems
\begin{theorem}
Given any  $R > 0$, the equation (\ref{E2})  has a solution   \ $(\lambda,u)$, \ where \ $\lambda \ \in 
\ \mathbb{R}$, \ $u \ \geq \ 0$, \ $u \ \not\equiv \ 0$ \ is a \ $C^{1,\alpha}$ function, \ $0 < \ \alpha \ < \  1$,
\ verifying \ $B(u) \ = \ R$ \ and \ $I(u)\ = \ \mu_R$, \ if \ $f \in C^{\infty}$ \ satisfies the following conditions: \\
$(p_1)$ \ $f(t,x)$\ is a stricly increasing odd function on \ $t$; \\
$(p_2)$ \ There exist constants \ $b \ >\ 0$\ and\ $0\ <\ \rho  \ <\ p^{\ast} - 1$ \ such that
 \ $|f(t,x)|\ \leq \ b\left(1 + |t|^{\rho}\right)$.
\end{theorem}

\begin{theorem}
The equation \ (\ref{E2}) \ has a solution \ $(\lambda,u), \ \lambda \in \mathbb{R},\ u \ \in \ C^{1,\alpha}(M)$ \ for some 
\ $0 \ < \ \alpha \ < \ 1$, \ $u \ \geq \ 0$ \ and \ $u \ \not\equiv \ 0$ \ if \ $f(t,x)$ \ satisfies to the following 
properties: \\
$(p_1)$ \ $f(t,x)$ is a stricly increasing odd function on \ $t$; \\
$(p_3)$ \ There exist positive constants  \ $b$ \ and \ $c$ \ such that $ |f(t,x)| \ \leq \ b \ + \ c|t|^{p^{\ast} - 1};$
\\
$(p_4)$ \ $\ddd\lim_{t \to 0^{+}}\inf \frac{1}{t^{p^{\ast} - 1}} [\inf_{x \in M}f(t,x)] \ = \ + \ \infty.$ 
\vspace{0.3cm}

\noindent  The  function \ $u$ is  strictly positive and increasing for \ $\lambda \ \geq \ 0$.
\end{theorem}

We list the article by O. Druet \cite {druet4}, where he studied a generalization of $ (Y) $ for a more general 
operator (the $ p $-Laplacian), as the article by Aubin \cite {aubin3}, to obtain a solution \ $ ( \lambda, u), \ \lambda \in 
\mathbb {R} \, \mbox { and} \, u \in H^p_1 $,  to the equation (\ref{E2}).

\noindent To find such a solution used as a main tool, the  Lagrange Multipliers's Theorem, which can be used because 
of the nature of the equation.

\section{Generalization of a nonlinear differential equation}

In this section we will work with a generalization of paper of Aubin \cite{aubin3}, where he has considered the differential  
equation  (\ref{E1}), namely $\Delta u \ +\  a(x)u \ = \ \lambda f(u,x),$
 on a compact Riemannian manifold \ $(M,g)$, where \ $a(x)$\ is a function \ $C^{\infty}$ \ on \ $M$\ and \ $f(u,x)$\ is a  \ $C^{\infty}$ function
 \ on \ $\mathbb{R}\times M$. 

In his paper, Aubin showed that, under certain conditions on \ $f(u,x)$, \ the equation \ (\ref{E1}) \ has a
 regular solution whenever \ $f(u,x)$ \  satisfies the increasement condition: there are two positive constants
  \ $b$ \ \ and \ $\rho$ \ such that 
$|f(t,x)| \ \leq \ b\left(1 + |t|^{\rho}\right),$ 

\noindent where \ $0\  < \ \rho \ \leq \ (n + 2)/(n - 2) \ = \ 2^{\ast}  -  1$, \ $2^{\ast}  =  (2n)/(n - 2)$. 

We will use the method in \cite{aubin3} to generalize the below equation, in the sense of that the operator will be the $p$-Lapacian.
 For this,
  by the lack of compactness of Sobolev embedded for the critical case (Theorem of compact embedded of Kondrakov)
we split the development into two cases: subcritical case
($ 0 < \rho < p^{\ast} - 1$) \  \ and the critical case \ ($\rho  =  p^{\ast}  -  1). $\ This kind of equation was studied by
  many authors in the  Euclidean context. In the Riemannian context we refer mainly to the Druet's article
\cite{druet4} which we extracted regularities' theorems and Maximun principles were used.

Let \ $(M,g)$ \ be a compact Riemannian manifold, \ $n$-dimensional,\ $n \ \geq \ 3$ and \ $p \ \in \ (1,n)$. 

We are interested in the following generalization of the equation (\ref{E1}):

\noindent  We look for solutions \ $u \ \in \ H^p_1(M)\cap C^0(M)$\  and\ $\lambda \ \in \ \mathbb{R}$\ for the equation (\ref{E2}), 
namely
\vspace{-0.5mm}
$$\Delta_p u \ + \ a(x)u^{p -1} \ = \ \lambda f(u,x) $$ 

\noindent  where \ $|f(t,x)| \ \leq \ b\left(1 + |t|^{\rho}\right)$ ,
 \ $0 < \ \rho \ \leq \ p^{\ast} -1$,\ \ $p^{\ast} = pn/(n - p)$ 

\noindent and \ $\Delta_p u \ = \ - \ div\left(|\nabla u|^{p - 2}\nabla u\right)$ \ is the \ $p$-Laplacian \ of \ $u$. 

{\bf Remark.} {\it If $p = 2$, the equation (\ref{E2}) became to (\ref{E1}), since $\Delta_2 u = \Delta u$.}

\subsection{Subcritical case}


In this section we will study the equation \ (\ref{E2}) \ in the subcritical case, i. e., where 
$ 0 <\rho <p^{\ast} - 1 $. \ The goal is to obtain a solution as the limit of a minimizing sequence for
the invariant \ $ \mu_R $  that, after using the Dominated Convergence Theorem of Lebesgue,  can be directly used
in the subcritical case because of the  compact embedded of Sobolev, in this case, the convergence to a solution follows
easily from  Lagrange Multipliers's Theorem.

For the proof of Theorem 1, we need the following lemma:

\begin{lemma}
If \ $f(t,x)$ \ satisfies the condition \ $(p_1)$, \ then

\noindent  $(i)$ \ $F(t,x)$ \ is a non negative  and \ $C^{\infty}$ function.

\noindent  $(ii) \ F(0,x) \ = \ 0 \ \mbox{and} \ F(\infty,x) \ = \ \infty.$

\noindent  $(iii) \ F(t,x)$ \ is an increasing function for \ $t \geq 0$.

\noindent  $(iv) \ F(t,x) \ = \ F(|t|,x) \ \forall \ t.$
\end{lemma}

\noindent  {\bf Proof of Lemma 1:} 

\noindent  $(i)$ \ As $F(t,x) \ = \ \ddd\int^t_0{f(s,x)ds}$ and \ $f$ \ is of $C^{\infty}$ class,
 we have that \ $F \in C^{\infty}$.

As \ $f$ \ is increasing and odd, \ $f(0,x) \ = \ 0$ \ and if \ $t \geq 0, \ \ F(t,x) \ \geq \ 0.$

Now, if \ $t < 0$, \ take \ $m > 0$ \ such that \ $t = - m$. \ So
\vspace{-0.5mm}
\begin{eqnarray*}
F(t,x) & = & \ddd\int^t_0{f(s,x)ds} = \ddd\int^{- m}_0{f(s,x)ds} \\
       & = & - \ddd\int^0_{- m}{f(s,x)ds} = \ddd\int^0_{- m}{f(- s,x)ds}, \ \mbox{taking} \ z = - s, \\
       & = & - \ddd\int^0_m{f(z,x)dz} = \ddd\int^m_0{f(z,x)dz} \geq 0.
\end{eqnarray*}
\fimaf

\noindent  $(ii) \ F(0,x) \ = 0$ \ is concluded directly by definition.

Taking \ $t_1 > 0$, \ then
\vspace{-0.5mm}
\begin{eqnarray*}
F(\infty,x) & = & \ddd\int^{\infty}_0{f(s,x)ds} = \ddd\int^{t_1}_0{f(s,x)ds} + \ddd\int^{\infty}_{t_1}{f(s,x)ds} \\
            & \geq & A + f(t_1,x)\ddd\int^{\infty}_{t_1}{ds} = \infty, 
\end{eqnarray*}

\noindent where \ $A = \ddd\int^{t_1}_0{f(s,x)ds}.$ \fimaf

\noindent  $(iii)$ \ If \ $0 \leq t_1 < t_2$, \ then

$F(t_2,x) = \ddd\int^{t_2}_0{f(s,x)ds} = \ddd\int^{t_1}_0{f(s,x)ds} + \ddd\int^{t_2}_{t_1}{f(s,x)ds} > 
\ddd\int^{t_1}_0{f(s,x)ds} = F(t_1,x).$ \fimaf

\noindent  $(iv)$ \ If \ $t \geq 0, \ F(t,x) = F(|t|,x).$

\noindent  If \ $t < 0$,
\vspace{-0.5mm}
\begin{eqnarray*} 
F(t,x) & = & \ddd\int^t_0{f(s,x)ds}, \ \mbox{taking} \ s = - z, \\
       & = & - \ddd\int^{- t}_0{f(- z,x)dz} = \ddd\int^{- t}_0{f(z,x)dz} = F(|t|,x).
\end{eqnarray*}       
\fim

\noindent  {\bf Proof of Theorem 1:}

By using item \ $(iv)$ \ of Lemma 1, we can consider  
$\mu_R = \ddd\inf_{u \in \mathcal{H}_R}{I(u)},$ 

\noindent  where 
$\mathcal{H}_R  = \{u \in H^p_1(M); \ u \geq 0 \ \mbox{and} \ B(u) = R \}.$

{\bf Remark}
By items \ $(ii) \ \mbox{and} \ (iii)$ from Lemma 1, clearly \ $\mathcal{H}_R \neq \emptyset$.

The proof of the theorem follows in several steps:

{\bf Claim 1}
There exist \ $N > 0$ \ such that, if \ $u \in \mathcal{H}_R$, \ then \ $\|u\|_1 \leq N$.

Firstly fix a \ $t_o > 0$. Then \ $\forall \ u \ \in \ \mathcal{H}_R$
\vspace{-0.5mm}
\begin{eqnarray*}
\|u\|_1 = \ddd\int_M{u dV} & = & \ddd\int_{\{u < t_o\}}{u dV} + \ddd\int_{\{u \geq t_o\}}{u dV} \\
   & \leq & t_o vol(M) + \ddd\int_{\{u \geq t_o\}}{u dV}.
\end{eqnarray*}   

\noindent  For $u \ \geq \ t_o \ > \ 0$, we have \ $f(u,x)\ \geq \ f(t_o,x) \ \geq \ \eta \
 = \ \ddd\inf_{x \in M}{f(t_o,x)} \ > \ 0$ \ by \ $(p_1)$. Whence 
\vspace{-0.8cm}

\begin{eqnarray*}
R = B(u) & = & \ddd\int_M{F(u,x)dV} \\
  & \geq & \ddd\int_{\{u \geq t_o \}}{F(u,x)dV} \\
  & = & \ddd\int_{\{u \geq t_o \}}{\left[\ddd\int_{t_o}^{u(x)}{f(t,x)dt}\right]dV} \\ 
  & \geq & \ddd\int_{\{u \geq t_o \}}{\left[\ddd\int_{t_o}^{u(x)}{f(t_o,x)dt}\right]dV} \\
  & \geq & \ddd\int_{\{u \geq t_o \}}{\left[\ddd\int_{t_o}^{u(x)}{\eta dt}\right]dV} \\ 
  & = & \eta \ddd\int_{\{u \geq t_o \}}{(u(x) - t_o)dV} \\
  & = & \eta \ddd\int_{\{u \geq t_o \}}{u(x)dV} \ - \ \eta t_ovol(\{u \geq t_o \}) \\
  & \geq & \eta \ddd\int_{\{u \geq t_o \}}{u(x)dV} \ -\  \eta t_ovol(M).
\end{eqnarray*}

So, $\ddd\int_{\{u \geq t_o \}}{u(x)dV} \ \leq \ \ddd\frac{R}{\eta} \ + \ t_o vol(M)$, 
\vsp

\noindent  where $\{u \geq t_o \} \ = \ \{x \ \in \ M \ ; \ u(x) \ \geq \ t_o \}$ \ and
 \ $vol(X)$ \ is the volume of $X \ \subseteq \ M$. 
\vsp

\noindent  Then, 
\vspace{-1cm}

\begin{eqnarray*}
 \|u \|_1  =  \ddd\int_M{u(x)dV} & = & \ddd\int_{\{u \geq t_o \}}{u(x)dV} \ + \ \ddd\int_{\{u < t_o\}}{u(x)dV} \\
 & \leq & \ddd\frac{R}{\eta} \ + \ 2t_ovol(M) \ = \ N. 
\end{eqnarray*}
\fimaf

{\bf Claim 2}
$\mu$ is finite.

Indeed, by using the  below inequality  (see \cite{aubin6}), for every \ $\epsilon \ > \ 0$ \ corresponds a \ $C(\epsilon) \ > 0$
such that
\begin{eqnarray}
\label{2.1}
\ddd\int_M{|u|^pdV} \ \leq \ \epsilon \ddd\int_M{|\nabla u|^pdV} \ + \ C(\epsilon)\left[\ddd\int_M{|u|dV}\right]^p 
\ \  \forall \ u \ \in \ H^p_1.
\end{eqnarray}

Therefore, for \ $u \in \mathcal{H}_R$, \ we have 
\vspace{-0.5mm}
$$I(u) \ = \ \ddd\int_M{|\nabla u|^pdV} \ + \ \ddd\int_M{a|u|^pdV} \ \geq \ \ddd\int_M{|\nabla u|^pdV} \ + 
\ \inf_Ma\ddd\int_M{|u|^pdV}.$$

If \ $\inf_Ma \ \geq \ 0$, \ we have \ $I(u) \ \geq \ 0$ \ and, consequently, \ $\mu \ \geq \ 0$.

If $\inf_Ma \ < \ 0$, \ by using \ (\ref{2.1}) \ and Claim 1, we have that
\vspace{-0.8cm}

\begin{eqnarray*}
I(u) & \geq & \ddd\int_M{|\nabla u|^pdV} \ + \ \epsilon \inf_Ma \ddd\int_M{|\nabla u|^pdV} \ + 
\ (\inf_Ma)C(\epsilon)\left[\ddd\int_M{|u|dV}\right]^p \\
& \geq & (1 \ + \ \epsilon \inf_Ma)\ddd\int_M{|\nabla u|^pdV} \ + \ (\inf_Ma)C(\epsilon)N^p \\
 & \geq & (\inf_Ma)C(\epsilon)N^p \ > \ - \infty, 
\end{eqnarray*}

\noindent  Since \ $\epsilon \ > \ 0$ \ is such that \ $1 \ + \ \epsilon \ddd\inf_Ma \ > \ 0$. \ What conclude the Claim 2.
 \fimaf 
\vsp

Consider now a sequence \ $(u_j) \ \in \ H^p_1$\ , \ $u_j \ \geq \ 0$ \ , \ $B(u_j) \ = \ R$ \ and
 \ $I(u_j) \ \rightarrow \ \mu_R $ \ when \ $j \ \rightarrow \ \infty$ \ (minimizing sequence).

{\bf Claim 3}
$(u_j)$ \ is bounded in \ $H^p_1$.

Indeed, as \ $I(u_j) \ \rightarrow \ \mu_R$, \ there exist \ $K \ > \ 0$ \ 
such that \ $|I(u_j)| \ \leq \ K \ \forall \ j$. Then by \ (\ref{2.1}) and Claim 1, respectively 
\vspace{-0.8cm}

\begin{eqnarray*}
\|\nabla u_j\|^p_p & = & I(u_j) \ - \ \ddd\int_M{a|u_j|^pdV} \ \leq \ K \ + \ \sup_M|a|\ddd\int_M{|u_j|^pdV}\\
& \leq  & K \ + \ \epsilon\sup_M|a|\|\nabla u_j\|^p_p \ + \ C(\epsilon)\sup_M|a|\|u_j\|^p_1 \\
& \leq & K \ + \ \epsilon\sup_M|a|\|\nabla u_j\|^p_p \ + \ C(\epsilon)\sup_M|a|N^p. 
\end{eqnarray*}

\noindent  So 
\vspace{-0.5mm}
$$(1 \ - \ \epsilon \ddd\sup_M|a|)\|\nabla u_j\|^p_p \ \leq \ K \ + \  C(\epsilon)\ddd\sup_M|a|N^p.$$ 

\noindent  Then, taking \ $\epsilon \ > \ 0$ such that \ $1 \ - \ \epsilon \ddd\sup_M|a| \ > \ 0$, \ we obtain
\begin{eqnarray}
\label{2.2}
\|\nabla u_j\|^p_p \ \leq \ C, 
\end{eqnarray}

\noindent  where \ $C \ > \ 0$ is a positive  constant.
\vsp

Therefore, by \ (\ref{2.1}), (\ref{2.2}) \ and Claim 1, we conclude the proof of Claim 3. \fimaf
\vsp

Now, as \ $H^p_1$ is reflexive and the Sobolev's embedded \ $H^q_1 \ \hookrightarrow \ L^s$  is compact for 
\ $1 \ \leq \ s \ < \ p^{\ast}$, the Claim 3 guarantees the existence of a subsequence \ $(u_i)$ of \ $(u_j)$ and
 \ $u \ \in \ H^p_1$ such that
\vspace{-0.5mm}
$$u_i \ \rightharpoonup \ u  \ \ \mbox{in} \ \ H^p_1, \eqno{(A_1)}$$
$$u_i \ \longrightarrow  \ u \ \mbox{in} \ \ L^s, \ 1 \ \leq \ s \ < \ p^{\ast} \ \ \mbox{and} \eqno{(A_2)}$$
$$u_i \ \longrightarrow \ u \ \ \mbox{a.e. in} \ M. \eqno{(A_3)}$$

By \ $(A_1)$ \ and \ $(A_2)$ \ \ $I(u) \ \leq \ \ddd\lim_{i \to \infty}\inf I(u_i) \ = \ \mu_R$. 
\vsp

By \ $(A_3)$ \ $u \ \geq \ 0$ \ a.e. in \ $M$. \ From \ $(A_2)$ and \ $(p_2)$ 
we can use the  Lebesgue Dominated Convergence's  Theorem   (see \cite{aubin6}) to conclude that \ $B(u) \ = \ R$.
\vsp

Hence \ $I(u) \ = \ \mu_R$, \ $u \ \geq \ 0$ \ com \ $u \ \not\equiv \ 0$.
\vsp

So, as \ $B$ \ and \ $I$ \ $\in \ C^1(H^p_1)$, \ taking \ $S \ = \ \{v \ \in \ H^p_1 \ ; \ B(v) \ = \ R \} $, 
\ we have that \ $B^{'}(v) \ \neq \ 0$ \ for every \ $v \ \in \ S$ \ and \ $u \ \in \ S$ \ is such that \ $I(u) \ = 
\ \ddd\inf_{v \in S}I(v)$. \ Then, by Lagrange Multipliers's Theorem  (see \cite{aubin6}),
 exist \ $\xi \ \in \ \mathbb{R}$ \ such that $I^{'}(u) \ = \ \xi B^{'}(u) $ namely 
$$p\ddd\int_M{|\nabla u|^{p - 2}\nabla u \nabla \varphi dV} \ + \ p\ddd\int_M{u^{p - 1}\varphi dV}
 \ = \ \xi \ddd\int_M{f(u,x)\varphi dV} \ \ \forall \ \varphi \ \in \ H^p_1.$$ 

In other words, \ $u$ \ is a solution of the equation $\Delta_p u \ + \ a(x)u^{p -1} \ = \ \lambda f(u,x)$,
in the weak sense, where \ $\lambda \ = \ \xi/p$. 

Finally, by \ $(p_2)$  we can use the Regularity Theorem (see  \cite{druet4})
 to conclude that exist \ $0 \ < \alpha \ < \ 1$ such that \ $u \ \in \ C^{1,\alpha}(M)$. \fim

{\bf Remark}
If \ $\lambda \ \geq \ 0$, \ by the Strong maximum principle's Theorem and (see \cite{druet4}) \ $u \ > \ 0$ \ in \ $M$.

\section{Critical case}
\mbox{}

We will study now the equation  (\ref{E2}), where  \ $ \rho  =  p^{\ast}  -  1$.  The problem here is the lack of compactness for
Sobolev's embedded  when \ $s \ = \ p^{\ast}$ \ (Kondrakov's theorem of embbed) and, to circumvent this difficulty,
it will be added an additional condition on \ $f(u,x)$. The goal is bring down the critical level of   $f $ and use Theorem 1. 

\noindent  {\bf Proof of Theorem 2:}

For each \ $ m \in \mathbb{N}^{\ast}$, \ define
\vspace{-0.5mm}
$$f_m(t,x) \ = \ signal(t).|f(t,x)|^{m/(m+1)}.$$

\noindent  Then, \ $f_m(t,x)$ \ is an odd function and strictly increasing in \ $t$ \ and, by \ $(p_3)$, 
\ it satisfies \ $(p_2)$ \ of Theorem 1.

Fixing \ $R \ > \ 0$ \ (to be clarified further on), as \ $f(t,x)$ \ satisfies \ $(p_1)$, by items
 \ $(ii) \ \mbox{and} \ (iii)$ of Lemma 1, exist \ $\nu \ \in \ \mathbb{R}$, \ $\nu \ > \ 0$ such that
$$\ddd\int_M{F(\nu,x)dV} \ = \ R \ , \ \ \ \mbox{where} \ \ F(\nu,x) \ = \ \ddd\int_0^{\nu}{f(t,x)dt}.$$

Now define
$$F_m(t,x) \ = \ \ddd\int^{t}_0{f_m(s,x)ds}$$

\noindent and
$$B_m(u) \ = \ \ddd\int_M{F_m(u(x),x)dV}.$$

Putting
$$R_m \ = \ \ddd\int_M{F_m(\nu,x)dV},$$

$$\mathcal{H}_m  = \{u \in H^p_1(M); \ u \geq 0 \ \mbox{and} \ B_m(u) = R_m \}$$

\noindent  and
$$\mu_m \ = \ \ddd\inf_{u \in \mathcal{H}_m}{I(u)},$$

\noindent  then, by Theorem 1, for each \ $m \ \in \ \mathbb{N}^{\ast}$,
 \ exist a function \ $u_m \ \in \ C^{1,\alpha}$, \ $u_m \ \geq \ 0$, \ $u_m \ \not\equiv \ 0$ 
\ and a real number \ $\lambda_m$ \ satisfying 
\vspace{-0.5mm}
\begin{eqnarray}
\label{2.3}
\Delta_pu_m \ + \ a(u_m)^{p-1} \ = \ \lambda_m |f(u_m,x)|^{m/(m+1)}, 
\end{eqnarray}

\noindent  because \ $signal(u_m) = 1.$ Moreover, \ $u_m$ performs
$$B_m(u_m) \ = \ \ddd\int_M{F_m(u_m(x),x)dV} \ = \ R_m$$
\noindent  and
$$\mu_m \ = \ I(u_m).$$

{\bf Claim 4}
$(u_m)$ \ is bounded in \ $H^p_1$.

Indeed, as
$$F_m(\nu,x) \ = \ \ddd\int_0^{\nu}{|f(t,x)|^{m/(m+1)}dt} \ \leq \ \nu \ + \ F(\nu,x)$$
we have
\begin{eqnarray}
\label{2.4}
R_m \ \leq \ \nu .vol(M) \ + \ R  \ \ \forall \ m. 
\end{eqnarray}
On the other hand, fixing \ $t_o > 0 \ \mbox{and} \ \eta > 0$ \ like in proof of Claim 1
\vspace{-0.5mm}
\begin{eqnarray*}
\|u_m\|_1 = \ddd\int_M{u_m dV} & = & \ddd\int_{\{u_m < t_o\}}{u_m dV} + \ddd\int_{\{u_m \geq t_o\}}{u_m dV} \\
   & \leq & t_o vol(M) + \ddd\int_{\{u_m \geq t_o\}}{u_m dV}.
\end{eqnarray*}   

\noindent  For $u_m \ \geq \ t_o \ > \ 0$, \ we have \ $f(u_m,x)\ \geq \ f(t_o,x) \ \geq \ \eta \ > \ 0$. \ Whence 
$$|f(t_o,x)|^{m/(m + 1)} \ \geq \ \eta^{m/(m + 1)}$$

\noindent  and 
\begin{eqnarray*}
R_m = B_m(u_m) & = & \ddd\int_M{F_m(u_m,x)dV} \\
  & \geq & \ddd\int_{\{u_m \geq t_o \}}{F_m(u_m,x)dV} \\
  & = & \ddd\int_{\{u_m \geq t_o \}}{\left[\ddd\int_{t_o}^{u_m(x)}{|f(t,x)|^{m/(m + 1)}dt}\right]dV} \\ 
  & \geq & \ddd\int_{\{u_m \geq t_o \}}{\left[\ddd\int_{t_o}^{u_m(x)}{|f(t_o,x)|^{m/(m + 1)}dt}\right]dV} \\
  & \geq & \ddd\int_{\{u_m \geq t_o \}}{\left[\ddd\int_{t_o}^{u_m(x)}{\eta^{m/(m + 1)} dt}\right]dV} \\ 
  & = & \eta^{m/(m + 1)} \ddd\int_{\{u_m \geq t_o \}}{(u_m(x) - t_o)dV} \\
  & = & \eta^{m/(m + 1)} \ddd\int_{\{u_m \geq t_o \}}{u_m(x)dV} \ - \ \eta^{m/(m + 1)} t_o vol(\{u_m \geq t_o \}) \\
  & \geq & \eta^{m/(m + 1)} \ddd\int_{\{u_m \geq t_o \}}{u_m(x)dV} \ -\  \eta^{m/(m + 1)} t_o vol(M).
\end{eqnarray*}

Thus, $\ddd\int_{\{u_m \geq t_o \}}{u_m(x)dV} \ \leq \ R_m\eta^{-m/(m + 1)} \ + \ t_o vol(M)$ 
\vsp

\noindent and, by \ (\ref{2.4}), \ exist a \ $C > 0$ \ such that  
\vspace{-0.5mm}
\begin{eqnarray}
\label{2.5}
\|u_m\|_1 \ \leq \ C, \ \forall \ m. 
\end{eqnarray}

Now, as  
$$\mu_m \ = \ I(u_m) \ \leq \ I(\nu) \ = \ \nu^p\ddd\int_M{a(x)dV}$$

\noindent we obtain
\vspace{-0.5cm}

\begin{eqnarray*}
\|\nabla u_m\|^p_p & = & I(u_m) \ - \ \ddd\int_M{a|u_m|^pdV} \\
 & \leq & \nu^p\ddd\int_M{a(x)dV} \ + \ \ddd\sup_M|a|\int_M{|u_m|^pdV} \\
 & \leq & \nu^p\ddd\int_M{a(x)dV} \ + \ \epsilon \sup_M|a|\|\nabla u_m\|_p^p \ + \ C(\epsilon)\sup_M|a|\|u_m\|^p_1, 
\end{eqnarray*}

\noindent  where \ $\epsilon \ > \ 0$ \ and \ $C(\epsilon) \ > \ 0$ \ came from \ (\ref{2.1}).

By taking \ $\epsilon \ > \ 0$ \ small enough so that \ $1 \ - \ \epsilon \sup_M|a| \ > \ 0$ 
\ we have, by \ (\ref{2.5}), \ that exist \ $C \ > \ 0$ \ such that 
\vspace{-0.5mm}
\begin{eqnarray}
\label{2.6}
\|\nabla u_m\|^p_p \ \leq \ C \ \ \forall  \  m. 
\end{eqnarray}

Finally, by using \ (\ref{2.1}), \ (\ref{2.5}) \ and \ (\ref{2.6}) \ we conclude the Claim 4. \fimaf

{\bf Claim 5}
$(\lambda_m)$ \ is bounded in \ $\mathbb{R}$.

Indeed, multiplying \ (\ref{2.3}) \ by\ $u_m$ \ and integrating on \ $M$, \ we obtain
\vspace{-0.5mm}
\begin{eqnarray}
\label{2.7}
I(u_m) \ = \ \lambda_m \ddd\int_M{u_m|f(u_m,x)|^{m/(m+1)}dV}. 
\end{eqnarray}

On the other hand, as \ $\|u_m\|_{H^p_1} \ \leq \ C$, \ there is \ $A \ > \ 0$ such that 
\vspace{-0.5mm}
\begin{eqnarray}
\label{2.8}
|I(u_m)| \ \leq \ A \ \ \forall \ m. 
\end{eqnarray}

By \ $(p_1)$ we have
\vspace{-0.5mm}
\begin{eqnarray*}
R_m & = & \ddd\int_M{F_m(u_m,x)dV}  =  \ddd\int_M{\left[\int_0^{u_m}{|f(t,x)|^{m/(m+1)}dt}\right]dV} \\
 & \leq & \ddd\int_M{\left[\int_0^{u_m}{|f(u_m,x)|^{m/(m+1)}dt}\right]dV} \\
 & = & \ddd\int_M{u_m|f(u_m,x)|^{m/(m+1)}dV}. 
\end{eqnarray*}

Now, by using \ (\ref{2.7}), \ (\ref{2.8}) \ and the above expression, we obtain
$$A \ \geq \ |I(u_m)| \ = \ |\lambda_m|\ddd\int_M{u_m|f(u_m,x)|^{m/(m+1)}dV} \ \geq \ |\lambda_m |R_m.$$
Namely,
\vspace{-0.5mm}
\begin{eqnarray}
\label{2.9}
|\lambda_m |R_m \ \leq \ A \ \ \forall \ m. 
\end{eqnarray}

Furthermore, when \ $m \ \rightarrow \ \infty$, \ $f^{m/(m+1)}(t,x) \ \rightarrow \ f(t,x)$ and
the convergence is dominated by \ $1 \ + \ f(t,x),$ it is integrable over \ $[0,\nu]$, whence \ $F_m(\nu,x) 
\ \rightarrow \ F(\nu,x)$. \ And the convergence is dominated by \ $\nu \ + \ F(\nu,x)$. \ Then we have (see \cite{aubin6}) 
\vspace{-0.2cm}
$$R_m \ \rightarrow \ R \ \ \mbox{when} \ m \ \rightarrow \ \infty.$$
As\ $R \ > \ 0$, \ we  can assume that there is \ $C_o \ > \ 0$ \ such that \ $R_m \ > \ C_o$ \ \ $\forall \ m$.
\vsp

So, by  (\ref{2.9}), \ $|\lambda_m| \ \leq \ \ddd\frac{A}{C_o}$ \ , \ give us the proof of Claim 5. \fimaf
\vsp

As \ $H^p_1$ is reflexive  the Sobolev's embedded 
\ $H^p_1 \ \hookrightarrow \ L^s$ is compact for \ $1 \ \leq \ s \ < p^{\ast}$, 
\ from Claims 4 and 5, there are \ $(u_m)$  subsequence of \ $(u_m)$, \ $(\lambda_m)$ 
subsequence of \ $(\lambda_m)$, \ $u \ \in \ H^p_1$ \ and \ $\lambda \ \in \ \mathbb{R}$ such that 
\vspace{-0.5mm}
$$u_m \ \rightharpoonup \ u \ \ \mbox{in} \ H^p_1, \eqno{(B_1)}$$
$$u_m \ \longrightarrow \ u \ \mbox{in} \ L^s, \ \ 1  \leq s <  p^{\ast}, \eqno{(B_2)}$$
$$u_m \ \longrightarrow \ u \ \ \mbox{a.e. in} \ M \ \mbox{and} \eqno{(B_3)}$$
$$\lambda_m \ \longrightarrow \ \lambda. \eqno{(B_4)}$$

{\bf Remark}
We are using in the proofs the same notation to denote a subsequence.

With this, \ $u \ \geq \ 0$ \ and \ $|f(u_m,x)|^{m/(m+1)} \ \longrightarrow \ f(u,x)$ \ a.e. in \ $M$.

{\bf Claim 6}
$|f(u_m,x)|^{m/(m+1)}$ is bounded in \ $L^{p^{\ast}/(p^{\ast}-1)}$.

Indeed, by Hölder's inequality 
\vspace{-0.5mm}
\begin{eqnarray*}
\ddd\| |f(u_m,x)|^{m/(m+1)}\|^{(m+1)/m}_{p^{\ast}/(p^{\ast}-1)} & = & \ddd\| f(u_m,x)\|_{[m/(m+1)][p^{\ast}/(p^{\ast}-1)]} \\
 & \leq & vol(M)^{(p^{\ast} - 1)/(m + 1)p^{\ast}} \ddd\| f(u_m,x)\|_{p^{\ast}/(p^{\ast}-1)} \\
 & \leq & C \ddd\| f(u_m,x)\|_{p^{\ast}/(p^{\ast}-1)}
\end{eqnarray*} 

\noindent  and by \ $(p_3)$
\vspace{-1cm}

\begin{eqnarray*}
\ddd\| f(u_m,x)\|_{p^{\ast}/(p^{\ast}-1)} & = & \left[\ddd\int_M{|f(u_m,x)|^{p^{\ast}/(p^{\ast}-1)}dV}\right]^{(p^{\ast}-1)/p^{\ast}} \\
 & \leq & \left[\ddd\int_M{(b_1 \ + \ c_1|u_m|^{p^{\ast}})dV}\right]^{(p^{\ast}-1)/p^{\ast}} \\
 & \leq & C \ + \ C\|u_m\|^{p^{\ast}-1}_{p^{\ast}} \ \leq \ C
\end{eqnarray*}

\noindent  this last inequality is due to Claim 4 and  $H^p_1
  \hookrightarrow  L^{p^{\ast}}$. \ $b_1$, \ $c_1$ \ are  positive constants and \ $C$ \ represent several
 positive constants, not necessarily the same.

We conclude the proof of Claim 6. \fimaf
\vsp

Consequently, (see \cite{aubin4}),  considering a subsequence, 
\vspace{-0.5mm}
\begin{eqnarray}
\label{2.10}
|f(u_m,x)|^{m/(m+1)} \ \rightharpoonup \ f(u,x) \ \mbox{em} \ L^{p^{\ast}/(p^{\ast}-1)}. 
\end{eqnarray}

Analogously, by Claim 4, \ $|\nabla u_m|^{p-2}\nabla u_m$ is bounded in \ $L^{p/(p-1)}$. Then, considering a subsequence
$|\nabla u_m|^{p-2}\nabla u_m \ \rightharpoonup \ \Sigma  \ \mbox{in} \ L^{p/(p-1)},$ for some \ $\Sigma \ \in \ L^{p/(p-1)}$.
\vsp

Now, by using \ (\ref{2.3}), \ $(p_3)$, \ $(B_2)$ \ and \ $(B_4)$ we conclude that \ $div(|\nabla u_m|^{p-2}\nabla u_m)$
 is bounded in \ $L^1$, we have that  \ $\Sigma \ = \ |\nabla u|^{p-2}\nabla u$ (see \cite{druet4}). Therefore,
\vspace{-0.5mm}
\begin{eqnarray}
\label{2.11}
|\nabla u_m|^{p-2}\nabla u_m \ \rightharpoonup \ |\nabla u|^{p-2}\nabla u \ \ \mbox{in} \ L^{p/(p-1)}. 
\end{eqnarray}

To conclude the proof of the Theorem 2, we remember from (\ref{2.3}) that
\vspace{-0.5mm}
$$\ddd\int_M{|\nabla u_m|^{p-2}\nabla u_m\nabla \varphi dV} \ + \ \ddd\int_M{a(u_m)^{p-1}\varphi dV} \ 
= \ \lambda_m \ddd\int_M{|f(u_m,x)|^{m/(m+1)}\varphi dV}, \  \forall \ \varphi \ H^p_1.$$

Taking \ $m \ \rightarrow \ \infty$, and using \ $(B_2)$, $(B_4)$, \ (\ref{2.10}) \ and \ (\ref{2.11}), we obtain 
\vspace{-0.5mm}
$$\ddd\int_M{|\nabla u|^{p-2}\nabla u\nabla \varphi dV} \ + \ \ddd\int_M{au^{p-1}\varphi dV} \ =
 \ \lambda \ddd\int_M{f(u,x)\varphi dV}, \ \ \ \ \forall \ \varphi \ \in \ H^p_1.$$

Namely, \ $u$ \ is a solution (in the weak sense) of the equation (\ref{E2}).

To regularize the solution we use the hypothesis\ $(p_3)$ (see \cite{druet4}). With this, there is some
 \ $0 \ < \ \alpha \ < \ 1$ \ such that \ $u \ \in \ C^{1,\alpha}(M)$.

As we already know that \ $u \ \geq \ 0$, \ to finish the proof of the theorem we have to show that \ $u \ \not\equiv \ 0$.
\vsp

By \ $(B_1)$ and \ $(B_2)$, we have that
\vspace{-0.5mm}
\begin{eqnarray}
\label{2.12}
I(u) \ \leq \ \lim_{m \to \infty}\inf I(u_m). 
\end{eqnarray}

For some function \ $u_o \ \in \ H^p_1$, \ $u_o \ \geq \ 0$, \ $u_o \ \not\equiv \ 0$, if we have 
\ $I(u_o) \ \leq \ 0$, \ then for each \ $m$, \ there is \ $k_m \ > \ 0$ \ such that
 \ $B(k_mu_o) \ = \ R_m$ \ and\ $I(k_mu_o) \ = \ (k_m)^pI(u_o) \ \leq \ 0$ \ (see Lemma 1). Then
 \ $\mu_m \ = \ I(u_m) \ \leq \ 0$ forall\ $m \ \geq \ 1$ \ and, using \ (\ref{2.12}), \ $I(u) \ \leq \ 0$. 

If \ $I(u) \ = \ 0$ \ we conclude that  $\mu_m \ = \ 0, \ \lambda \ = \ \lambda_m \ = \ 0 \ \mbox{and} \ u_m \ \not\equiv \ 0$ 
 satisfies 
$\Delta_pu_m \ + \ a(x)(u_m)^{p-1} \ = \ 0,$
  this concludes the proof of theorem.

But, if \ $I(u) \ < \ 0$, \ we have that \ $u \ \not\equiv \ 0$, this also prove the theorem.
\vsp

Let us prove, then, the case where \ $I(u_m) \ > \ 0$ \ for all \ $m \ \geq \ 1$.

By \ $(p_3)$, \ we have
\vspace{-0.5mm}
$$|f(t,x)|^{m/(m+1)} \ \leq \ b_1 \ + \ c_1 |t|^{[m/(m+1)][p^{\ast}-1]} \ \leq \ b_1 \ + \ c_1 \ + \ c_1|t|^{p^{\ast}-1}$$

\noindent  where \ $b_1$ \ and \ $c_1$ are  positive constants. Thus, considering \ $b_2 \ = \ b_1 \ + \ c_1$, \ we obtain
\vspace{-0.5mm}
\begin{eqnarray}
\label{2.13}
R_m \ = \ \ddd\int_M{\left[\ddd\int_0^{u_m}{|f(t,x)|^{m/(m+1)}dt}\right]dV} \ \leq \ b_2\|u_m\|_1 \ +
 \ \ddd\frac{c_1}{p^{\ast}}\|u_m\|^{p^{\ast}}_{p^{\ast}}. 
\end{eqnarray}

As \ $H^p_1 \ \hookrightarrow \ L^{p^{\ast}}$, \ there is \ $K \ \mbox{and} \ D \ > \ 0$ \ such that 
\vspace{-0.5mm}
$$\|\varphi \|^p_{p^{\ast}} \leq \ K\|\nabla \varphi \|^p_p \ + \ D\|\varphi \|^p_p \ \ \forall \ \varphi \ \in \ H^p_1.$$

\noindent  From this fact
\vspace{-0.5mm}
\begin{eqnarray}
\label{2.14}
\|\varphi \|^{p^{\ast}}_{p^{\ast}} \leq \ \left[\ K\|\nabla \varphi \|^p_p \ + 
\ D\|\varphi \|^p_p \ \right]^{p^{\ast}/p} \ \  \forall \ \varphi \ \in \ H^p_1. 
\end{eqnarray}

Then, by \ (\ref{2.13}) and \ (\ref{2.14}) we have
\vspace{-0.5mm}
\begin{eqnarray}
\label{2.15}
R_m \ - \ b_2\|u_m\|_1 \ \leq \ \ddd\frac{c_1}{p^{\ast}}\left[\ K\|\nabla u_m \|^p_p \ + \ D\|u_m \|^p_p \ \right]^{p^{\ast}/p}. 
\end{eqnarray}

If \ $R_m \ - \ b_2\|u_m\|_1 \ < \ 0$, then \ $\|u_m\|_1 \ > \ \ddd\frac{R_m}{b_2}$, what give us, by \ $(B_2)$ \ and by
 \ $R_m \ \rightarrow \ R \ > \ 0$, \ that \ $\|u\|_1 \ \geq 
\ \ddd\frac{R}{b_2} \ > \ 0$, \ in other words, \ $u \ \not\equiv \ 0$.
\vsp

Now, if \ $R_m \ - \ b_2\|u_m\|_1 \ \geq \ 0$, we have \ $1  \ - \ \ddd\frac{b_2}{R_m}\|u_m\|_1 
\ \geq \ 0$ \ and by \ (\ref{2.15}) \ we obtain
\vspace{-0.5cm}

\begin{eqnarray}
\label{2.16}
\left(\ddd\frac{R_mp^{\ast}}{c_1}\right)^{p/p^{\ast}}\left(1 \ - \ \ddd\frac{b_2}{R_m}\|u_m\|_1\right) 
& \leq & \left(\ddd\frac{R_mp^{\ast}}{c_1}\right)^{p/p^{\ast}}\left(1 \ - \ \ddd\frac{b_2}{R_m}\|u_m\|_1\right)^{p/p^{\ast}}
 \nonumber \\
 & \leq & K\left[ \ I(u_m) \ - \ \ddd\int_M{a|u_m|^pdV} \ \right] \ + \ D\|u_m\|^p_p \nonumber \\
 & \leq & \mu_m K \ + \ D_o\|u_m\|^p_p
\end{eqnarray}

\noindent  where \ $D_o \ > \ 0$.
 
{\bf Claim 7}
There is \ $\epsilon \ > \ 0$ \ such that forall \ $m \ \geq \ 1$ \ and a convenient \ $R \ > \ 0$  
$$K\mu_m \ < \ \left(\ddd\frac{Rp^{\ast}}{c_1}\right)^{p/p^{\ast}} \ - \ 2\epsilon.$$

Indeed, by \ $(p_4)$ \ there is a sequence of real numbers \ $\nu_i \ > \ 0$ \ such that \ $\nu_i
 \ \rightarrow \ 0$, \ when \ $i \ \rightarrow \ \infty$, \ and \ $f(t,x) \ > \ it^{p^{\ast}-1}$ \ forall
 \ $t \ \in \ (0,\nu_i)$. \ This implies that 
$$F(\nu_i,x) \ = \ \ddd\int_0^{\nu_i}{f(t,x)dt} \ > \ \ddd\frac{i}{p^{\ast}}(\nu_i)^{p^{\ast}}$$
and, consequently,
$$R_i \ = \ \ddd\int_M{F(\nu_i,x)dV} \ > \ \ddd\frac{i}{p^{\ast}}(\nu_i)^{p^{\ast}}vol(M).$$

Taking
$$(R_m)_i \ = \ \ddd\int_M{F_m(\nu_i,x)dV},$$ 
\vspace{-0.5mm}
$$(\mathcal{H}_m)_i \ = \ \{u \in  H^p_1(M);\ u \geq 0 \ \mbox{and} \ B_m(u)\ = \ (R_m)_i \}$$

\noindent  and 
\vspace{-0.5mm}
$$(\mu_m)_i \ = \ \ddd\inf_{u \in (\mathcal{H}_m)_i}{I(u)},$$

\noindent  we obtain
\vspace{-0.5mm}
$$(\mu_m)_i \ \leq \ I(\nu_i) \ = \ (\nu_i)^p\ddd\int_M{a(x)dV}.$$ 

With this
\vspace{-0.5mm} 
$$ \ddd\frac{(\mu_m)_i}{(R_i)^{p/p^{\ast}}} \ \leq \ \left[(\nu_i)^p\int_M{a(x)dV}\right] 
/\left[\left(\frac{i}{p^{\ast}}\right)^{p/p^{\ast}}(\nu_i)^p.vol(M)^{p/p^{\ast}}\right] \ \longrightarrow
 \ 0 \ \mbox{when} \ \ i \ \rightarrow \ \infty .$$

{\bf Remark}
Remember that \ $R_i \ > \ 0$ \ and we are considering the case where \ $(\mu_m)_i \ > \ 0$ \ forall \ $m$ \ and  $i \ \geq \ 1$.

Hence, 
$$\ddd\frac{K(\mu_m)_i}{\left[(R_ip^{\ast})/c_1\right]^{p/p^{\ast}}} \ \longrightarrow \ 0 \ \ \mbox{when} \ \ i 
\ \rightarrow \ \infty \ , \  \ \forall \ \ m \ \geq \ 1.$$

Then, for a big enough \ $i$, taking \ $R \ = \ R_i$ \ and
 \ $\mu_m \ = \ (\mu_m)_i$, \ we have that there is \ $\epsilon_o \ > \ 0$ \ such that, forall \ $m \ \geq \ 1$
\vspace{-0.3cm}
$$\ddd\frac{K(\mu_m)}{\left[(Rp^{\ast})/c_1\right]^{p/p^{\ast}}} \ < \ 1 \ - \ \epsilon_o$$

\noindent  and taking \ $2\epsilon \ = \ \epsilon_o\left[(Rp^{\ast})/c_1\right]^{p/p^{\ast}}$ \ 
we conclude proof of Claim 7. \fimaf
\vsp

Now, by Claim 7 and the fact that \ $R_m \ \rightarrow \ R$ \ when \ $m \ \rightarrow \ \infty$, 
\ after some \ $m_o$
\vspace{-0.3cm}
$$K\mu_m \ + \ \epsilon \ < \ \left(\ddd\frac{R_mp^{\ast}}{c_1}\right)^{p/p^{\ast}}$$

\noindent  and, by using \ (\ref{2.16}), we obtain 
\vspace{-0.3cm}

$$\left(K\mu_m \ + \ \epsilon \right)\left(1 \ - \ \ddd\frac{b_2}{R_m}\|u_m\|_1\right) \ \leq \ K\mu_m \ + \ D_o\|u_m\|^p_p.$$

Then, 
\vspace{-0.3cm}

$$\epsilon \ - \ \left(K\mu_m \ + \ \epsilon \right)\ddd\frac{b_2}{R_m}\|u_m\|_1 \ \leq \ D_o\|u_m\|^p_p \ $$

\noindent  consequently
\vspace{-0.5mm}
\begin{eqnarray}
\label{2.17}
\epsilon \  \leq \ \left(K\mu_m \ + \ \epsilon \right)\ddd\frac{b_2}{R_m}\|u_m\|_1 \ + \  D_o\|u_m\|^p_p 
\end{eqnarray}

\noindent  and since that \ $\mu_m \ > \ 0$, \ $R_m \ \rightarrow \ R$ \ when \ $m \ \rightarrow \ \infty$, 
\ by \ $(B_2)$ \ and \ (\ref{2.17}) \ $u \ \not\equiv \ 0$.
\vsp

Finally, if \ $\lambda \ \geq \ 0$, \ by Strong  maximum principle (see \cite{druet4}), \ $u \ > \ 0$. \fim

{Silva, C.R}
{{INSTITUTO DE CI\^ENCIAS EXATAS E DA TERRA\\
 CAMPUS UNIVERSIT\'ARIO DO ARAGUAIA\\UNIVERSIDADE FEDERAL DE MATO GROSSO
 \\  78698-000, Pontal do Araguaia,  MT, Brazil}  }

{carlosro@ufmt.br}


{Romildo Pina}
{
Instituto de Matem\'atica e Estat\'istica\\
Universidade Federal de Goi\'as - UFG\\
74690-900, Goi\^ania, GO, Brazil}

{romildo@ufg.br}
\thanks{}

{Marcelo Souza}
{
Instituto de Matem\'atica e Estat\'istica\\
Universidade Federal de Goi\'as - UFG\\
74690-900, Goi\^ania, GO, Brazil}

{msouza@ufg.br}
\thanks{}


\begin{thebibliography}{9}

\markright{\footnotesize{References}}

\bibitem{aubin0} AUBIN, T. {\it Métriques riemanniennes et courbure}. Journal of Differential Geometry, vol. 4, n. 4, p. 383-424, 1970.

\bibitem{aubin3} --------------. {\it Equations différentielles non linéaires}. Bulletin des Sciences Mathématiques, vol. 99, 
p. 201-210, 1975.

\bibitem{aubin2} --------------. {\it Equations différentielles non linéaires et problème de Yamabe concernant la 
courbure scalaire}. Journal de Mathématiques Pure et Appliquées, vol. 55, p. 269-296, 1976.

\bibitem{aubin4} --------------. {\it Problème isopérimétriques et espaces de Sobolev}. Journal of Differential 
Geometry, vol. 11, p. 573-598, 1976.

\bibitem{aubin6} --------------. {\it Some nonlinear problems in riemannian geometry}. Springer Monographs in Mathematics, 1998.

\bibitem{azorero1} AZORERO, J. G., \& ALONSO, I. P. {\it Existence and nonuniqueness for the p-Laplacian: nonlinear eigenvalues}. Communications in Partial Differential Equations, vol. 12, p. 1389-1430, 1987.

\bibitem{brezis1} BRÉZIS, H., \& NIRENBERG, L. {\it Positive solutions of nonlinear elliptic equations involving critical Sobolev exponents}. Communications on Pure and Applied Mathematics, vol. 36, p. 437-477, 1983.


\bibitem{correa} F.J. CORR\^EA,  J.V. GONCALVES \& A.L. MELO,  {\it On positive radial solutions of quasilinear elliptic equations}. Nonlinear Analysis, vol. 52, p. 681-701, 2003.

\bibitem{demegel} DEMEGEL, F., \& HEBEY, E. {\it On some nonlinear equations involving the p-Laplacian with critical Sobolev growth}. Advances in Differential Equations, vol. 3, n. 4, p. 533-574, 1998.

\bibitem{djadli} DJADLI, Z. {\it Nonlinear elliptic equations with critical Sobolev exponent on compact riemannian manifolds}. Calculus of Variations and Partial Differential Equations, vol. 8, p. 293-326, 1999.


\bibitem{druet4} DRUET, O. {\it Generalized scalar curvature type equations on compact riemannian manifolds}. Proceedings of the Royal Society of Edinburgh, vol. 130 A, p. 767-788, 2000.

\bibitem{olimpio} MIYAGAKI, O. H. {\it On a class of semilinear elliptic problems in $\mathbb{R}^n$ with critical growth}. Nonlinear Analysis, Theory, Methods and Applications, vol. 29, n. 7, p. 773-781, 1997.

\bibitem{schoen} SCHOEN, R. {\it Conformal deformation of a riemannian metric to constant scalar curvature}. Journal of Differential Geometry, vol. 20, p. 479-495, 1984.

\bibitem{silva}  SILVA, C. R. {\it Algumas equações diferenciais não-lineares em variedades riemannianas compactas, UnB. thesis (2004)}


\bibitem{silva2} SILVA, C. R. {\it On the study of Existence of solutions for a
class of equations with critical Sobolev exponent on compact Riemannian Manifold. Mat. Contempor\^anea (2014), vol. {\bf  43}, p. 223-246}

\bibitem{trudinger} TRÜDINGER, N. S. {\it Remarks concerning the conformal deformation of riemannian structures on compact manifolds}. Ann. Scuola Norm. Sup. Pisa, vol. 3, n. 22, p. 265-274, 1968.

\bibitem{yamabe} YAMABE, H. {\it On a deformation of riemannian structures on compact manifolds}. Osaka Math. J., vol. 12, p. 21-37, 1960.

\end{thebibliography}
\end{document}